\documentclass[12pt,a4paper]{amsart}
\usepackage{amscd,amssymb,amsopn,amsmath,amsthm,graphics,amsfonts,enumerate,verbatim,calc}

\headsep=1cm \numberwithin{equation}{section}
\newtheorem{theorem}{Theorem}[section]
\newtheorem{proposition}[theorem]{Proposition} 

\newtheorem{corollary}[theorem]{Corollary}
\newtheorem{remark}[theorem]{Remark}
\newtheorem{example}[theorem]{Example}
\newtheorem{definition}[theorem]{Definition}

\numberwithin{equation}{section}

\def\<{\left < }
\def\>{\right >}
\def\({\left ( }
\def\){\right )}
\def\e{\eqref}

\def\2{\#_2 }

\begin{document}

\title{Geometry and topology of maximal antipodal sets and related topics}
\author{Bang-Yen Chen}

\address{Department of Mathematics, Michigan State University, 
East Lansing, MI 48824, USA}
 \email{chenb@msu.edu}

\begin{abstract} Maximal antipodal sets of Riemannian manifolds were introduced by the author and T. Nagano in  [Un invariant g\'eom\'etrique riemannien, C. R. Acad. Sci. Paris S\'er. I Math.  295 (1982), no. 5, 389--391]. Since then maximal antipodal sets have been studied by many mathematicians and they shown that maximal antipodal sets are related to several important areas in mathematics. 
The main purpose of this paper is thus to present a  comprehensive  survey on geometry and topology of maximal antipodal sets and also on their applications to several related topics. 
 \vspace{0.2 cm} 
   \end{abstract}

\maketitle

\section{Maximal antipodal sets and two-numbers}

The notion of  maximal antipodal sets of Riemannian manifolds was introduced  in \cite{CN82}. The primeval concept of  maximal antipodal sets is the notion of antipodal points on a circle. 
For a circle $S^1$ in the Euclidean  plane $\mathbb E^2$, the {\it antipodal point} $q$ of a point $p\in S^1$ is the point on $S^1$ which is diametrically opposite to $p$. 

In Riemannian geometry, a geodesic in a Riemannian manifold $M$ is a curve which yields locally the shortest distance between any two nearby points.
Since a closed geodesic in a Riemannian manifold $M$ is isometric to a planar circle,  antipodal points  can be defined for every closed geodesic in $M$, i.e., a point $q$  in a closed geodesic is called an {\it antipodal point} of another point $p$ on the closed geodesic if the distance $d(p,q)$ between $p$ and $q$ on the two arcs connecting $p$ and $q$ are equal. 
For simplicity, a closed geodesic in a Riemannian manifold is also called a {\it circle} in this article.

A subset $S$ of a Riemannian manifold $M$ is said to be an {\it antipodal set} if any two points in $S$ are antipodal on some circle of $M$. An antipodal set in a connected Riemannian manifold $M$ is called a {\it maximal antipodal set} if it doesn't lie in any antipodal set  as a proper subset. The  supremum of the cardinality of all maximal antipodal set of $M$ is called the {\it two-number} of $M$, denote by $\#_2M$. If an antipodal set $S$ of $M$ satisfies
$\# S=\#_{2}M$, then $S$ is called a {\it great antipodal set} or a {\it 2-set}.

\vskip.05in
If a Riemannian manifold $M$ contains no closed geodesics, we put  $\#_2M=0$. 
On the other hand, for any compact Riemannian manifold $M$ we have 
\begin{align}\label{1.1} \#_2M\geq 2,\end{align} 
since every compact connected Riemannian manifold contains at least one close geodesic (see \cite{LF51}). Clearly, inequality \e{1.1} is sharp, because $\#_2 S^{n}=2$ for every standard $n$-sphere $S^{n}$. A clear proof of the finiteness of the 2-number for compact connect Riemannian manifolds was given  by M. S. Tanaka and H. Tasaki \cite{TT13}.

\section{$(M_{+},M_{-})$-theory for compact symmetric spaces}\label{S2}

\subsection{E. Cartan's classification of irreducible compact symmetric spaces}\label{S2.1}

In 1926,  \`Elie Cartan achieved his classification of symmetric spaces by  reducing the problem to the classification of simple Lie algebras over  real field ${\mathbb R}$, a problem which Cartan himself solved earlier in 1914.

\vskip.03in
\`E. Cartan proved that simply-connected irreducible symmetric  spaces of compact  type consist of the following four families:

\begin{itemize}
\item[{\rm (1)}] Classical simple Lie groups  $SO(n), SU(n), Sp(n)$. 

\vskip.04in
\item[{\rm (2)}] Exceptional simple Lie groups   $E_{6}, E_{7}, E_{8}, F_{4},G_{2}$. 

\vskip.04in
\item[{\rm (3)}] Eight classes of classical symmetric spaces $AI(n), AII(n),AIII(p,q), BDI(p,q),BDII(n)$,

$DIII(n),CI(p),CII(p,q)$ corresponding to the classical groups $SO(n), SU(n), Sp(n)$.

\vskip.04in
\item[{\rm (4)}] Twelve  exceptional symmetric spaces $EI,EII,EIII,EIV,EV,EVI,EVII,EVIII, EIX$,

$FI,FII$ and $GI$ corresponding to the exceptional simple groups $E_{6}, E_{7}, E_{8}, F_{4},G_{2}$. 
 \end{itemize}
 
\vskip.04in
 Since the discovery  by \'E Cartan, symmetric spaces, a distinguished
class of Riemannian manifolds, attracted the attention of numerous mathematicians from various fields such as differential geometry, algebraic topology, representation theory and harmonic analysis.

 \subsection{$(M_{+},M_{-})$-theory}\label{S2.2}

A geometric new approach to compact symmetric spaces, the $(M_{+},M_{-})$-theory, was developed by the author and T. Nagano during the 1970--1980s (see \cite{C87,CLN80,CN77,CN78,CN82,CN88}). This theory was also known today as the Chen--Nagano theory  in some literatures (see, e.g., \cite{Ta14,Ta15}).

The fundamental principle of this theory is that a pair of a polar
and the corresponding meridian determines a compact Riemannian symmetric space. One of advantages of this theory is that it is useful for inductive arguments on polars or meridians.

\vskip.05in

Assume that $o$ is a point of a compact symmetric space, say $M=G/K$. A connected component of the fixed point set $F(s_{o},M)\smallsetminus \{o\}$ of the symmetry $s_{o}$ at $o$ is called a {\it polar} of $o$.  We denote it by $M_{+}$ or $M_{+}(p)$ if $M_{+}$ contains a point $p$.  
 
The following  propositions from \cite{CN78} is quite useful (also \cite[page 15]{C87}).

\begin{proposition}\label{P:2.1} Let $M=G/K$ be a compact symmetric space. Then, for each antipodal point $p$ of $o\in M$, the isotropy subgroup $K$ at $o$ acts transitively on the polar $M_+(p)$. Further, we have $K(p)=M_+(p)$ and $K(p)$ is connected. Hence,  $M_+(p)=K/K_p$, where $K_p=\{k\in K: k(p)=p\}$.
\end{proposition}
 
If a polar consists of a single point, then it is called a {\it pole}.

\vskip.05in
 
 \begin{proposition}\label{P.2.2} Under the hypothesis of Proposition \ref{P:2.1}, the normal space to $M_+(p)$ at $p\in M$ is the tangent space of a connected complete totally geodesic submanifold  $M_-(p)$. Thus we have
  \begin{align}\dim M_{+}(p)+\dim M_{-}(p)=\dim M.\end{align}
\end{proposition}

We call $the M_{-}(p)$ the {\it meridian} of $o$ at $p$. For meridians, we have 

 \begin{proposition}\label{P:2.3} For each antipodal point $p$ of $o$ in a compact symmetric space $M$, we have
 \begin{itemize}
 \item[{\rm (1)}]  $rk (M_-(p))=rk (M)$ and
\item[{\rm (2)}] $M_-(p)$ is a connected component of the fixed point set $F(s_{p}\circ s_{o},M)$ of $s_{p}\circ s_{o}$ through $p$.
\end{itemize}\end{proposition}

For a compact symmetric space $M$,  polars and meridians of $M$ are totally geodesic submanifolds; in fact, they are  compact symmetric spaces as well. 
Polars and meridians have been determined for every compact connected irreducible Riemannian symmetric space (see  \cite{C87,CN88,Na88,Na92}). 

One of the most important properties of polars and meridians is that $M$ is determined globally by any pair of $(M_{+}(p),M_{-}(p))$.

 \vskip.05in
Besides polars and meridians of  symmetric spaces of compact type, there exist another important totally geodesic submanifolds called centrosomes which are defined as follows.

\begin{definition} {\rm Let $o$ be a point of a compact connected Riemannian symmetric space $M$. If  $p$ is a pole of $o\in M$, then the {\it centrosome}   $C(o,p)$ of $\{o,p\}$  is the set consisting of the midpoints of all geodesics in $M$ joining $o$ and $p$.
A connected component of a centrosome is called a {\it centriole}.}
\end{definition}

\begin{remark} {\rm P. Quast described in \cite{Qu13} all centrioles in irreducible simply-connected compact  symmetric spaces of compact type in terms of the root system of the ambient space. He also investigated geometric properties of centrioles in  \cite{Qu13}.}
\end{remark}
  
 A connected component of the centrosome $C(o,p)$ is a totally geodesic submanifold of $M$.   Centrosomes play some important roles in topology as well. For example, centrosomes were used by J. M. Burns  to compute homotopy of compact symmetric spaces in \cite{Bu92}. 
 
  \vskip.05in
  The following result from \cite{CN88} characterizes poles in compact symmetric spaces (see also \cite{C87.2}).

\begin{proposition} 
 The following six conditions are equivalent to each other for two distinct points $o,p$ of a connected compact symmetric space $M=G_{M}/K_{G}$.

\begin{itemize} 

  \vskip.02in
  \item[(i)] $p$ is a pole of $o\in M$; 
  
 \vskip.02in
\item[(ii)] $s_{p}=s_{o}$; 
  
  \vskip.02in
\item[(iii)] $\{p\}$  is a polar of $o\in M$; 
  
 \vskip.02in
\item[(iv)] there is a double covering totally geodesic immersion $\pi=\pi_{\{o,p\}} :M\to M''$ with $\pi(p)=\pi(o)$;

 \vskip.02in
\item[(v)] $p$  is a point in the orbit $F(\sigma,G_{M})(o)$  of the group $F(\sigma,G_{M})$ through $o$, where $\sigma = {\rm ad}(s_{o})$;

 \vskip.02in
\item[(vi)] the isotropy subgroup of $SG_{M}$ at $p$ is that, $SK_{G}$ $($of $SG_{M}$ at $o)$, where $SG_{M}$ is the group generated by $G_{M}$ and the symmetries; $SG_{M}/G_{M}$ is a group of order $\leq 2$.
\end{itemize}  
\end{proposition}
  
For compact symmetric spaces, the author and T. Nagano proved the following result.

\begin{proposition} \label{P:2.9} \cite{CN88} For any compact symmetric space $M$,  the two-number $\#_2M$  is equal to the maximal possible cardinality $\#(A_2M)$ of a subset $A_2M$ of $M$ such that the point symmetry $s_x$ fixes every point of $A_2M$ for every $x\in A_2M$. 
\end{proposition}

 Proposition \ref{P:2.9} can be regarded as an alternative definition of 2-number for  symmetric spaces of compact type.

\begin{remark} {\rm T. Nagano and M. Sumi  \cite{NM88} proved  that the root system $R(M_-)$ of a meridian $M_-\ne M$ is obtained from the Dynkin diagram of the root system $R(M)$ of the compact symmetric space  $M$. Furthermore, they have determined  in \cite{NM88} all maximal totally geodesic spheres in $SU(n)$ by means of the $(M_+,M_-)$-theory.}
\end{remark}

\section{Descriptions of great antipodal sets}\label{S3}

In 1988, maximal antipodal sets of  compact symmetric spaces were determined and used by the author and T. Nagano  in  \cite{CN88} for determining 2-numbers of irreducible compact symmetries and also for simple Lie groups. In addition, they explicitly described antipodal sets in many compact Riemannian symmetric spaces, but  did not mention maximal antipodal sets for oriented real Grassmannian manifolds in \cite{CN88}.

\vskip.05in
On the other hand, there are quite many works done in recent years related to great antipodal sets. In particular, many authors have provided detailed descriptions of great antipodal sets for many symmetric spaces of compact  type. In this section, we will present some of their works in this respect.

\vskip.1in
\subsection{Great antipodal sets of classical  Lie groups}\label{S3.1}

In \cite{TT17}, M. S. Tanaka and H. Tasaki provided an explicit classification of maximal antipodal subgroups of compact classical Lie groups and of their factors by cyclic central subgroups. Their constructions of these subgroups are based on $D[4]$, i.e., the dihedral group  of order 8 or the automorphism group of a square in the plane. Also, maximum cardinalities of arbitrary antipodal sets in these compact Lie groups and factors are also calculated by them.

\vskip.05in

In \cite{TTY19},  M. S. Takano, H. Tasaki and O. Yasukura provided the classification of maximal antipodal subgroups of compact exceptional Lie group $G_{2}$ and compact symmetric space  $GI=G_{2}/SO(4)$ an explicit description of them by regarding as the automorphism group of the octonions $\mathcal O$. They also presented in  \cite{TTY19} the classification of maximal antipodal sets of  $GI$. Furthermore, they pointed out a relation between maximal antipodal sets of and those of the
oriented Grassmannian manifold $\widetilde G_{3}(\mathbb R^{7})$ by using the identification of with the set of associative 3-dimensional subspaces in ${\rm Im}(\mathcal O)$, which is a
totally geodesic submanifold of $\widetilde G_{3}(\mathbb R^{7})$.

\vskip.1in
\subsection{Great antipodal sets of exceptional Lie groups and  exceptional symmetric spaces}\label{S3.2}

 M. S. Tanaka, H. Tasaki and O. Yasukura  \cite{TTY22} explicitly described maximal antipodal sets of compact symmetric spaces related to $G_{2}$ by realizing it as the automorphism group of the octonions ${\mathcal O}$. Applying their explicit descriptions, they provided a close relation between maximal antipodal sets of the associative Grassmannian of the octonions and the Fano plane.

\vskip.05in
In \cite{Sa22.2}, Y. Sasaki provided an explicit classification of congruent classes of maximal antipodal sets of $F_{4}$ by using Jordan algebra $H_{3}({\mathcal O})$. Moreover, he explicitly classified congruent classes of maximal antipodal sets of $FI$.

\vskip.05in
Y. Sasaki classified in \cite{Sa22.2}  congruent classes of maximal antipodal sets of the exceptional Lie group$E_{6}$ and compact symmetric spaces of type  $EI, EII,EIII, EIV$ related to $E_{6}$ . Moreover, he gave realizations of these compact symmetric spaces by using some subalgebras of the complex exceptional Jordan algebra. In this realization, he explicitly described congruent classes of maximal antipodal sets of these compact symmetric spaces.

\vskip.1in
\subsection{Great antipodal sets of symmetric spaces of compact  type}\label{S3.3}

In \cite{Ta13}, H. Tasaki  described antipodal sets in oriented real Grassmannian manifolds, $\widetilde G_{k}(\mathbb R^{n})$. For a set $X$, the sets $$P_{k}(X)=\{\alpha\subset X:\#\alpha =k\}$$ are defined. The sets $P_{k}(n)=P_{k}(1,\ldots,n)$  were used to classify  maximal antipodal sets of $\widetilde G_{k}(\mathbb R^{n})$ for each $k\leq 4$. Furthermore, some arguments using for $k\leq 4$ are generalized to construct some maximal antipodal subsets for higher $k$.   In another article  \cite{Ta15},  H. Tasaki shown that great antipodal sets of  $\widetilde G_{5}(\mathbb R^{n})$ are unique up to isometries of $\widetilde G_{5}(\mathbb R^{n})$  for $n\geq 87$.

In \cite{TT20}, M. S. Tanaka and H. Tasaki  classified and explicitly described maximal antipodal sets of some compact classical symmetric spaces and those of their quotient spaces by making use of suitable embeddings of these symmetric spaces into certain compact classical Lie groups.  They also provided the cardinalities of maximal antipodal sets and they determined the maximum of the cardinalities and maximal antipodal sets whose cardinalities attain the maximum.

\vskip.05in
Very recently, J. Yu provided in \cite{Yu23}  explicit classification of maximal antipodal sets in any irreducible compact symmetric space except for spin and half spin groups, and some quotient symmetric spaces associated to them.

\vskip.1in
\subsection{Expansion of antipodal sets and homogeneous antipodal sets}\label{S3.4}

Y. Sasaki \cite{Sa21} introduced the notion of
connectedness of antipodal sets. Using connectedness, he defined a subgroup $G_{W}$ of the isometry group of a compact symmetric space $M$. He also constructed a method to build a bigger antipodal set from a given antipodal set via the subgroup $G_{W}$. 

\vskip.05in
An antipodal set $A\subset M$ is called {\it homogeneous} if there exists a subgroup of
the isometry group of $M$ acting on $A$ transitively.
In  \cite{Sa21}, Y. Sasaki proved that the connectedness is a sufficient condition that a maximal antipodal set is homogeneous.

\begin{remark} {\rm For further results on maximal antipodal sets of compact symmetric spaces, we refer to \cite{Bey18,LD14,Sa21,Sa22,Tas14,Tas17}.}
\end{remark}

\section{Links between two-numbers and topology}\label{S4}

Two-numbers link closely with topology. In this section, we provide several results in this respect.

\vskip.1in
\subsection{Two-numbers and Euler numbers}\label{S4.1}

In \cite{CN88}, the author and T. Nagano proved the following very simple link between 2-number and Euler number.

\begin{theorem} \label{T:4.1}  For  any compact symmetric space $M$, we have
\begin{align}\label{4.1}\#_2M\geq\chi(M),  \end{align}
where $\chi(M)$ denotes the Euler number of $M$.
\end{theorem}

The proof of this theorem was based on the $(M_+,M_-)$-theory in conjunction with a result of H. Hopf on fixed point sets and a result of Hopf and H. Samelson in \cite{HS}.

\vskip.05in
For any compact Hermitian symmetric space of semisimple type, the author and Nagano  proved the following.

\begin{theorem}\label{T:4.2} \cite{CN88}  For  any compact hermitian symmetric space $M$ of semi-simple type, we have
\begin{align}\label{4.2} \#_2M=\chi(M)=1+\sum \2 M_+.\end{align}
\end{theorem}
The proof of this theorem based heavily on the $(M_+,M_{-})$-theory as well as the Lefschetz
fixed point theorem in the version of M. F. Atiyah and I. M. Singer \cite{AS68}. 

\vskip.1in
The following result is an immediate consequence of Theorem \ref{T:4.2}.

\begin{corollary} For every complete totally geodesic hermitian subspace $B$ of a semi-simple hermitian symmetric space $M$, we have
$\chi(M)\geq\chi(B).$
\end{corollary}

\vskip.1in
\subsection{Two-numbers and covering maps}\label{S4.2}

 The author and T. Nagano discovered  in \cite{CN88} the following links between 2-numbers and covering maps for compact symmetric spaces.
   
\vskip.1in
For double coverings we have 
  
\begin{theorem} \label{T:4.4} If $M$ is a double covering of $M''$,  then $\#_{2}M\leq 2\#_2(M'').$
 \end{theorem}
 
\begin{remark} {\rm The inequality in Theorem \ref{T:4.4} is sharp, because the equality case holds for the group manifold $M=SO(2m)$ with $m>2$.}
\end{remark}

 For $k$-fold coverings with odd $k$, we have
 
\begin{theorem}  \label{T:4.6}  Let $\phi :\widetilde M\to M$ is a $k$-fold covering  between two compact symmetric spaces.  Then $\#_{2}\widetilde M=\#_2(M)$ whenever $k$ is odd.
 \end{theorem}

\vskip.1in
\subsection{A link between two-number and projective rank in algebraic geometry}\label{S4.3}

A. Fauntleroy  \cite{Fa93} defined  {\it projective rank}, denoted by $Pr(M)$, of a compact Hermitian symmetric space  $M$ as the maximal complex dimension of  totally geodesic complex projective spaces  of $M$.  
     
 \vskip.05in
C. U. S\'anchez and A. Guinta  \cite{SG02} proved  the following link between the two-number and the projective rank for irreducible Hermitian symmetric spaces of  compact type. 

 \begin{theorem} $Pr(M)\cdot rk(M)\leq \#_2(M)$ for any irreducible Hermitian symmetric space  of compact type.
  \end{theorem}

\subsection{Holomorphic two-number of compact hermitian symmetric space}\label{S4.4}

C. U. S\'anchez  defined in  \cite{Sa97}  the notion of {\it holomorphic two-number}, $\#_2^H(M)$, for a compact connected Hermitian manifold $M$ as the maximal possible cardinality of a subset $A_2$ such that for every pair of points $x$ and $y$ of $A_2$,  there exists a totally geodesic complex  curve of genus 0 in $M$ on which $x$ and $y$ are antipodal to each other.

 \vskip.05in
C. U. S\'anchez proved the following.

\begin{theorem} \label{T:4.8}  \cite{Sa97} $\#_{2}^H M=\#_2(M)\,$ for every compact hermitian symmetric space. 
\end{theorem}

By combining Theorem \ref{T:4.8} with our equality $\#_2M=\chi(M)$ from Theorem \ref{T:7.2},  one obtains  the following.
 
\begin{corollary}  \cite{Sa97}  $\#_2^H(M)$ $=\chi(M)$ for every compact hermitian symmetric space of semi-simple type.
\end{corollary}

\vskip.1in
\subsection{Symmetric $R$-spaces}\label{S4.5}

A  {\it R-space} (or a {\it real flag manifold}\/) is an orbit of the isotropy representation of a  symmetric space $G/K$ of compact type, where $G$ is a connected semisimple Lie group.
The notion of a symmetric $R$-spaces was introduced independent by T. Nagano \cite{Na65}  and M. Takeuchi \cite{Ta65} in 1965. By definition symmetric $R$-spaces are compact symmetric spaces which are at the same
time $R$-spaces. In fact, symmetric $R$-spaces admit a transitive action of a centre-free non-compact semisimple Lie group and the corresponding stabilizer of a point is a certain maximal parabolic subgroup. 

\vskip.05in
A  symmetric space $M$of compact type is said to have a {\it cubic lattice} if a maximal torus is isometric to the quotient of ${\mathbb E}^r$ by a lattice of ${\mathbb E}^r$ generated by an orthogonal basis of the same length. 
     O. Loos \cite{Lo85} gave another intrinsic characterization of symmetric $R$-spaces among all compact symmetric spaces with the property that the unit lattice of the maximal torus of the compact symmetric space (with respect to a canonical metric) is a cubic lattice. Loos' proof was based on the correspondence between the symmetric $R$-spaces and compact Jordan triple systems. 
     
\vskip.05in
S. Kobayashi and T. Nagano classified symmetric $R$-spaces in \cite{KNa64}. 
      The class of symmetric $R$-spaces consists of the following seven families:

\vskip.06in
\begin{itemize}
  \item[{\rm (a)}] All hermitian symmetric spaces of compact type

\vskip.02in
 \item[{\rm (b)}]  Grassmann manifolds $O(p+q)/O(p)\times O(q), Sp(p+q)/Sp(p)\times Sp(q)$

\vskip.02in
 \item[{\rm (c)}]  The classical groups $SO(m),\,U(m),\,Sp(m)$,

\vskip.02in
 \item[{\rm (d)}]  $U(2m)/Sp(m),\, U(m)/O(m)$,

\vskip.02in
 \item[{\rm (e)}]  $(SO(p+1)\times  SO(q+1))/S(O(p)\times O(q))$, where $S(O(p)\times O(q))$ is  the subgroup  of 
 
 $SO(p+1)\times SO(q+1)$ consisting of matrices of the form:
$$\begin{pmatrix} \epsilon& 0 & & \\ 0 & A&&\\ &&\epsilon& 0\\
&&0&B \end{pmatrix},\quad \epsilon=\pm1,\;\; A\in
O(p),\;\; B\in O(q),$$
 
\vskip.02in
 \item[{\rm (f)}]   Cayley projective plane $FII={\mathcal O}P^2$, and 

 \vskip.02in
\item[{\rm (g)}]  The three exceptional spaces $EIII=E_6/Spin(10)\times SO(2), EVII=E_7/E_6\times SO(2),$ 
 
 and $EIV=E_6/F_4.$
\end{itemize}

R. Bott  \cite{Bo59} used symmetric $R$-spaces to prove his famous periodicity theorem for the stable homotopy of classical Lie groups.   

Bott's original results may be succinctly summarized as 
 
 \begin{theorem}\label{T:4.11}  \cite{Bo70} The homotopy groups of the classical groups are periodic:
$$ \pi_i(U)=\pi_{i+2}(U),\; \pi_i(O)=\pi_{i+4}(Sp),\; \pi_i(Sp)=\pi_{i+4}(O)$$
for $i=0,1,\cdots$,
where $U$ is the direct limit defined by $U=\cup_{k=1}^\infty U(k)$ and similarly for $O$ and $Sp$.
\end{theorem}
 
\begin{remark} {\rm The second and third of these isomorphisms given in Theorem \ref{T:4.8} imply  the following 8-fold periodicity: }
$\pi_i(O)=\pi_{i+8}(O),\; \pi_i(Sp)=\pi_{i+8}(Sp),\;\;i=0,1,\cdots. $
\end{remark}

Also,  M. S. Tanaka and H. Tasaki proved in 2013  the following

 \begin{theorem}\label{TT13} \cite{TT13} Let $M$ be a symmetric $R$-space.
Then 
\begin{itemize}
  \item[{\rm (A)}] Every antipodal set is contained in a great antipodal set.
 \item[{\rm (B)}]Any two great antipodal sets are congruent,
where two subsets are congruent if they are 
  
  transformed to each other by
an element of the identity component of the isometry group.
  \end{itemize}\end{theorem}

\begin{remark} {\rm  Tanaka and Tasaki  \cite{TT13} also proved that there exists an  antipodal set of the adjoint group of $SU(4)$ which does not satisfy  condition (A) of Theorem \ref{TT13}.
Notice that the adjoint group of $SU(4)$  is a compact symmetric space, but not one of symmetric $R$-spaces.}
\end{remark}

\subsection{Intrinsically and extrinsically reflective submanifolds}\label{S4.6}

 A connected component $N$ of the fixed point set of an involutive isometry $\sigma$ of a Riemannian manifold $M$ is called a {\it reflective submanifold}. This isometry $\sigma$ is called the {\it reflection} of $M$ through $N$. Both polars and meridians of a compact symmetric spaces are reflective submanifolds.

\vskip.05in
A totally geodesic submanifold $M\subset \widetilde M$ of a submanifold $\widetilde M \subset \mathbb E^{m}$ of a Euclidean $m$-space $\mathbb E^{m}$
is called {\it extrinsically reflective}, if $M$ is a connected component of the intersection of $ \widetilde M$ with the fixed set of an involutive isometry of $\mathbb E^{m}$ that leaves $\widetilde M$ invariant.

A connected submanifold $P\subset \mathbb E^{m}$ of $\mathbb E^{m}$ is called an {\it extrinsically symmetric space} if for all $x\in P$ the submanifold $P$ is invariant under the reflections $\rho_{x}\in {\rm Isom}(\mathbb E^{m})$ through the affine normal space of $x+T^{\perp}_{x}(P)$ of $P$ at the point $x$, where $T^{\perp}_{x}(P)$  denotes the  normal space of $M$  in $\mathbb E^m$ at $x$.

\vskip.1in
In \cite{EQT15}, J.-H. Eschenburg, P. Quast and M. S. Tanaka proved the following. 

\begin{proposition} Every extrinsically reflective submanifold $M\subset P$ of an extrinsically symmetric space $P\subset \mathbb E^{m}$ is extrinsically symmetric in $\mathbb E^{m}$.\end{proposition}

\begin{proposition}  Any meridian $P_{-}$ of a compact extrinsically symmetric space $P\subset \mathbb E^{m}$ is itself extrinsically symmetric in $\mathbb E^{m}$.\end{proposition}

\begin{proposition} Every reflective submanifold of a compact extrinsically
symmetric space is actually extrinsically reflective, and thus extrinsically symmetric. In other
words, any reflective submanifold of a symmetric $R$-space is a symmetric $R$-space.
\end{proposition}

A connected submanifold $S$ of Riemannian manifold $M$ is called  (geodetically) {\it convex} if any shortest geodesic segment in $S$ is still shortest in $M$.

\vskip.05in
In 2012, P. Quast and M. S. Tanaka  proved the following.

\begin{theorem} \cite{QT12} Every reflective submanifold of a symmetric $R$-space is convex.\end{theorem}

\vskip.1in
\subsection{Links between two-numbers and homology}\label{S4.7}

Analogous to Theorem \ref{T:4.2}, M. Takeuchi proved in \cite{Ta89} the following. 

\begin{theorem}\label{T:4.13} \cite{Ta89}  For  any compact hermitian symmetric space $M$ of semi-simple type, we have
\begin{align} \#_2M=\chi(M)=1+\sum \2 M_+.\end{align}
\end{theorem}

 The $i$-th Betti number of  a manifold $M$  with coefficients in $\mathbb Z_2$ is the rank of the $i$-th homology group $H_{i}(M,\mathbb Z_2)$.
For any symmetric $R$-space, M. Takeuchi also proved  the following. 

\begin{theorem}\label{T:4.16} \cite{Ta89}  For  any symmetric $R$-space $M$, we have
\begin{align} \#_{2}M=\sum_{i\geq 0} b_i(M,\mathbb Z_2),  \end{align}
where $b_i (M,\mathbb Z_{2})$ is the i-th Betti number of $M$ with coefficients in $\mathbb Z_2$.
\end{theorem}
M. Takeuchi proved this theorem  by applying Theorem \ref{T:4.13} and a result of Chen--Nagano from \cite{CN88}  in conjunction with an earlier result of  Takeuchi in \cite{Ta65}.

\vskip.1in
\section{Antipodal sets and Borsuk-Ulam's theorem}\label{S5}

The following result in algebraic topology is well-known. 

\vskip0.1in
\noindent {\bf The Borsuk--Ulam Antipodal Theorem.} \cite{B33} {\it Every continuous function from an $n$-sphere $S^n$ into the Euclidean $k$-space $\mathbb E^k$ with $k\leq n$ maps some pair of antipodal points to the same point.}

\vskip0.1in
Obviously, Borsuk-Ulam's theorem fails for $k>n$, because $S^n$ can be embedded in $\mathbb E^{n+1}$. It is well-known that Borsuk-Ulam's theorem has numerous applications.
 For instance, H. Steinlein  provided in \cite{S85} a list of 457 publications involving various generalizations and/or applications of the Borsuk-Ulam theorem.

\vskip.05in
A  continuous function $f: M\to {\mathbb R}$ of a compact symmetric space $M=G/K$  is called {\it isotropic} if it is invariant under the action of the isotropic subgroup $K$. 

\vskip.05in
In  \cite{C17}, the author proved some Borsuk-Ulam's type theorems involving maximal antipodal sets of compact symmetric spaces as follows.

\begin{theorem}\label{T:5.1} Let  $f: M\to {\mathbb R}$ be an isotropic continuous function from a compact symmetric space. Then  $f$ maps  a great antipodal set of $M$ to the same point  in ${\mathbb R}$, whenever $M$ is one of the following spaces: Spheres; the projective spaces ${\mathbb F}P^n\, ({\mathbb F}={\mathbb  R}, {\mathbb C}, {\mathbb H})$;  the Cayley plane $FI\hskip-.02in I$; the exceptional spaces $EIV; EIV^*;GI$; and the exceptional Lie group $G_2$.
Hence,  $f: M\to {\mathbb R}$ maps a great antipodal set with $\#_2M$ elements to the same point  in ${\mathbb R}$.
\end{theorem}

The following example illustrates that the isotropic condition on $f$ in Theorem \ref{T:5.1} is necessary.

\begin{example} {\rm Let ${\mathbb R}P^2$ denote the real projective plane of curvature one. Then there exists a canonical double covering map $\pi:S^2(1) \to {\mathbb R}P^2.$ Assume that $S^2(1)$ is the unit sphere in $\mathbb E^3$ centered at the origin of $\mathbb E^3$. 
For each continuous function $f:{\mathbb R}P^2\to {\mathbb R}$, the lift $\hat f:S^2(1)\to {\mathbb R}$ of $f$ is an even function via the double covering $\pi$, so that $\hat f(-{\bf x})= \hat f({\bf x})$ for any ${\bf x}=(x,y,z)\in S^2(1).$

Conversely, for any continuous even function $h:S^2(1) \to \mathbb R$ of $S^2(1)$, $h$ induces a continuous function $\check{h}:{\mathbb R}P^2\to \mathbb R$ of ${\mathbb R}P^2$.
Let $h=(x-y)^2$. Then $h$ induces a function $\check{h}:{\mathbb R}P^2\to \mathbb R$ which does not map any  maximal antipodal set of ${\mathbb R}P^2$ to the same point in $\mathbb R$.}
\end{example}

Every compact symmetric spaces $M$ in the list of Theorem \ref{T:5.1} admits only a polar for  $o\in M$. Now, let us assume that $M$ is a compact symmetric space with multiple polars, said $M_+^1, M_+^2,\ldots,M_+^i$ for $o\in M$. In this case, let $\hat M_+$  denote a polar of $o\in M$ which has the maximal 2-number among all polars of $o$. For such compact symmetric spaces, we have the following result.

\begin{theorem} \cite{C17} Let  $f: M\to {\mathbb R}$ be an isotropic continuous function of a compact symmetric space $M$. If $M$ admits more than one polar,
then  $f$ maps an antipodal set of $M$ consisting of  $1+\#_2\hat M_+$ points of $M$ to the same point  in ${\mathbb R}$.
\end{theorem}

In particular, we have the following. 

\begin{theorem} \cite{C17} If  $f: E_8\to {\mathbb R}$ is  an isotropic continuous function, then  $f$ maps an antipodal set  of $E_8$ with  $392$ elements to the same point  in ${\mathbb R}$.
\end{theorem}

\begin{theorem} \cite{C17} If  $f: FI\to {\mathbb R}$ is  an isotropic continuous function of $FI$, then $f$ maps an antipodal set  of $FI$ with  $24$ elements to the same point  in ${\mathbb R}$.\end{theorem}

\section{Two-rank of Borel and Serre}\label{S6}

The 2-rank of a compact Lie group $G$ was introduced by A. Borel and J.-P. Serre \cite{BS53}.
 The 2-rank of  $G$, denoted by $r_2G$, is  the maximal possible rank of the elementary 2-subgroup of $G$. 

 \vskip.06in
Borel and Serre \cite{BS53} proved  the following:
\vskip.05in

\begin{itemize}
  \item[{\rm (a)}] 
 $rk(G)\leq r_2(G)\leq 2rk(G)$ and 
 
\vskip.05in
 \item[{\rm (b)}] $G$ has  2-torsion if $rk(G)<r_2(G)$,
\end{itemize}
where $rk(G)$ denotes the ordinary rank of $G$.
\vskip.05in

 In \cite{BS53}, Borel and Serre were able to determine  the 2-rank of  simply-connected simple Lie groups $SO(n), Sp(n), U(n), G_{2}$ and $F_{4}$. In addition,  they
proved that the exceptional Lie groups $G_{2}, F_{4}$ and $E_{8}$ have 2-torsion. On
the other hand, they  pointed out in  \cite[page 139]{BS53} that they were unable to determine the 2-rank for the exceptional simple Lie groups $E_{6}$ and $E_{7}$. 

\vskip.05in
After A. Borel and J.-P. Serre's paper,  2-ranks have been investigated by many mathematicians. For instance, it was shown that the 2-ranks have some links with commutative algebra.
Here, we  provide two of such links.

\vskip.1in
(i) Assume that $F$ is either a field or the rational integer ring $\mathbb Z$. Let 
$A=\sum_{i\geq 0} A_{i}$
 be a graded commutative $F$-algebra in sense of J. Milnor and J. Moore \cite{MM65}. If $A$ is connected, then it admits a unique augmentation $\varepsilon :A\to F$.   
 
 \vskip.06in
Put  $\bar A={\rm Ker}\, \varepsilon$. The $\bar A$ is called the {\it augmentation ideal} of $A$.
 A sequence of elements $\{x_{1},\ldots,x_{n}\in \bar A\}$ is said to be a {\it simple system of generators} if $\{x_{1}^{\epsilon_{1}}\cdots x_{n}^{\epsilon_{n}}: \epsilon_{i}=0 \hbox{ or }1\}$ is a module base of $A$.
For  a compact connected Lie group $G$, let us denote by $s(G)$ the number of generators of a simple system of the $\mathbb Z_{2}$-cohomology $H^{*}(G,{\mathbb Z}_{2})$ of  $G$. 

\vskip.05in
In \cite{Ko77}, A. Kono proved the following.

 \begin{theorem}\label{T:6.1}  \cite{Ko77} If $G$ is a connected compact Lie group, then
  the following three conditions are equivalent: 
 \begin{itemize}
  \item[{\rm (1)}] $s(G)\leq r_{2}G$; 
  
\vskip.02in
  \item[{\rm (2)}] $s(G)= r_{2}G$; 
  
 \vskip.02in
 \item[{\rm (3)}] $H^{*}(G,{\mathbb Z}_{2})$ is generated by universally transgressive elements. 
\end{itemize}
\end{theorem}

To prove Theorem \ref{T:6.1}, A. Kono applied P. May's spectral sequence \cite{Ma65},  S. Eilenberg and J. C. Moore's spectral sequence \cite{EM62} and also D. Quillen's result from \cite{Qu71.1}. 
\vskip.06in

In \cite{Ko77}, Kano also described some properties of compact Lie groups satisfying condition (3) in Theorem \ref{T:6.1} and provided some applications.

\vskip.1in

(ii)  The {\it Krull dimension} of a ring $R$ is the supremum of the number of strict inclusions in a chain of prime ideals, i.e., we say that a strict chain of inclusions of prime ideals of the form: 
 $${\mathfrak p}_{0}\subsetneq {\mathfrak p}_{1}\subsetneq \cdots \subsetneq {\mathfrak p}_{n}$$
 is of length $n$; i.e., it is counting the number of strict inclusions. 
   Given a prime ideal ${\mathfrak p}\subset R$, the {\it height} of ${\mathfrak p}$ is defined to be  the supremum of the set 
$$\{n\in {\mathbb N}: {\mathfrak p} \hbox{ is the supremum of a strict chain of length $n$}\},$$
and the Krull dimension is the supremum of the heights of all of its primes.
 
\vskip.06in
 Let $G$ be a compact Lie group. We put $H_{G}^{*}=H^{*}(BG; {\mathbb Z}_{2}),$ where $BG$ denotes the classifying space for $G$. Let $N_{G}^{*}\subset H_{G}^{*}$ be the ideal of nilpotent elements.  Then  $H_{G}^{*}/N_{G}^{*}=H_{G}^{\#}$  is a finitely generated commutative algebra.  
 
 \vskip.05in
In \cite{Qu71.1}, D. Quillen investigated the relationship between the finitely generated commutative algebra $H_{G}^{\#}$  and the  structure of the Lie group $G$. He proved that the Krull dimension of $H_{G}^{\#}$ is equal to the 2-rank of $G$ (under some suitable assumptions). 
 Quillen proved the result by calculating the mod 2 cohomology ring of extra special 2-groups. 
  Quillen's result gave rise to an affirmative answer to a conjecture of M. F.  Atiyah posed in \cite{At61}, and a conjecture of R. G. Swan given in  \cite{Sw71}.

 \section{Applications of maximal antipodal sets to Borel-Serre's problem}\label{S7}
 
 If $G$ is a connected compact Lie group, then by assigning $s_{x}(y)=xy^{-1}x$ to every point $x\in G$, we have $s^{2}_{x}=id_{G}$ to each point $x$. Thus, $G$ is a compact symmetric space with respect to a bi-invariant Riemannian metric. 
 
 \vskip.1in
\subsection{Links between two-numbers and 2-ranks}\label{S7.1}
  The author and Nagano proved the following  link between the 2-rank and the two-number of a connected compact Lie group.
 
 \begin{theorem}\label{T:7.1} \cite{CN82} Let $G$ be a connected compact Lie group. Then we have
 \begin{equation} \#_{2}G=2^{r_{2}G}.\end{equation}
 \end{theorem}
 
 For products of two compact Lie groups, we have the following result from \cite[Lemma 1.7]{CN88}.
 
\begin{theorem}\label{T:7.2} \cite{CN82} Let $G_{1}$ and $G_{2}$ be  connected compact Lie groups. Then 
 \begin{equation}\label{6.2} \#_{2}(G_{1}\times G_{2})=2^{r_{2}G_{1}+r_{2}G_{2}}.\end{equation}
\end{theorem}

\vskip.1in

 \subsection{2-ranks of classical groups}\label{S7.2} 
Applying Theorems \ref{T:7.1}, Theorem \ref{T:7.2} and  $(M_+,M_-)$-theory, the author and Nagano were able to determine the 2-ranks of all compact connected simple Lie groups in \cite{CN88}.  
Therefore, we have settled the problem of Borel-Serre for the determination of 2-ranks of all compact connected simple Lie groups.
 
\vskip.1in
 For classical groups we have:

\begin{theorem}\label{T:7.3}  Let $U(n)/{\mathbb Z}\mu$  be the quotient group of the unitary group $U(n)$ by the cyclic normal subgroup ${\mathbb Z}\mu$ of order $\mu$. Then we have 
  \begin{equation} r_{2}(U(n)/{\mathbb Z}\mu)= 
 \begin{cases} n+1 & \text{if $\mu$ is even and $n=2$ or $4$;}\\ n & \text{otherwise.}\end{cases} \end{equation}
 \end{theorem}

 \begin{theorem}\label{T:7.4}  For $SU(n)/{\mathbb Z}\mu$, we have 
 \begin{equation} r_{2}(SU(n)/{\mathbb Z}\mu)= 
 \begin{cases} n+1 & \text{for $(n,\mu)=(4,2)$;}
 \\ n  & \text{for $(n,\mu)=(2,2)$ or $(4,4)$;}
 \\ n-1 & \text{for the other cases.} \end{cases} \end{equation}
 \end{theorem}

 \begin{theorem}\label{T:7.5}  One has $r_{2}(SO(n))=n-1$ and, for $SO(n)^{*}$, we have
 \begin{equation} r_{2}(SO(n)^{*})= 
 \begin{cases} 4 & \text{for $n=4$;}\\ n-2 & \text{for $n$ even $>4$}.\end{cases} \end{equation}
 \end{theorem}

 \begin{theorem}\label{T:7.6}  Let $O(n)^{*}=O(n)/\{\pm 1\}$. We have
 
 {\rm (a)} $r_{2}(O(n))= n$;
 
 {\rm (b)} $r_{2}(O(n)^{*})$ is $n$ if $n$ is 2 or 4, while it is $n-1$ otherwise.
 \end{theorem}

 \begin{theorem}\label{T:7.7}  One has $r_{2}(Sp(n))=n$, and, for $Sp(n)^{*}$, we have
  \begin{equation} r_{2}(Sp(n)^{*})= \begin{cases} n+2 & \text{for $n=2$ or $4$}\\ n+1 & \text{otherwise}.
\end{cases} \end{equation}
Thus we also have
 \begin{equation}r_{2}(Sp(n)^{*})= r_{2}(U(n)/{\mathbb Z}_{2}) + 1 \end{equation} for every $n$. 
\end{theorem}

\vskip.1in

 \subsection{2-ranks of spinors,  semi-spinors and $Pin(n)$}\label{S7.3}
 
For $Spin(n)$ we have the following two results.

 \begin{theorem}\label{T:7.8}  We have
 \begin{equation}\notag
  r_{2}(Spin(n)) =\begin{cases} r+1 & \text{if  $\,n \equiv -1,0$ or {\rm 1 (mod 8)}}\\ 
 r &\text{otherwise},\end{cases}
 \end{equation} where $r$ is the rank of $Spin(n)$, $r = [\frac{n}{2}]$. 
 \end{theorem}

 \begin{theorem}\label{T:7.9} {\rm (PERIODICITY)} For $n\geq 0$, One has
$$ r_{2}(Spin(n+8))= r_{2}(Spin(n))+4$$
 
\end{theorem}

$Pin(n)$ was discovered  by M. F. Atiyah, R. Bott and A. Shapiro  while they studied Clifford modules in \cite{ABS64}. 
For the group $Pin(n)$ we have

 \begin{theorem}\label{T:7.10}  For $Pin(n)$ with $n\geq 0$, one has
 $$r_{2}(Pin(n))= r_{2}(Spin(n + 1)).$$  \end{theorem}
 
For the semi-spinor group $SO(4m)^{\#}=Spin(4m)/\{1, e_{((4m))}\}$, we have:

 \begin{theorem}\label{T:7.11}   We have
 \begin{equation}\notag
  r_{2}(SO(4m)^{\#}) =\begin{cases} 3 & \text{if  $m=1$}\\  6 &\text{if $m=2$,}\\ r+1 &\text{if $m$ is even $>2$, }\\
 r &\text{if $m$ is odd $>1$},\end{cases}
 \end{equation}
 where $r=2m$ is the rank of $SO(4m)^{\#}$.  
 \end{theorem}

\begin{remark}{\rm The 2-rank of $Spin(16)$ and of $SO(16)^{\#}$ were obtained  in \cite{Ad87} independently by J. F. Adams.  His method of proof was completely different from ours given in \cite{CN88}.}
\end{remark}

\vskip.1in

 \subsection{2-ranks of exceptional groups}\label{S7.4}

For exceptional Lie groups we have the following.

\begin{theorem}\label{T:7.13}   One has $r_{2}E^{*}_{6}=6$. \end{theorem}

\begin{theorem}\label{T:7.14} One has $$r_{2}G_{2}=3,\;\;  r_{2} F_{4}=5,\;\;  r_{2}E_{6}=6,\;\; r_{2}E_{7}=7
,\;\; r_{2}E_{8}=9$$ for  simply-connected exceptional simple Lie groups $G_{2}, F_{4}, E_{6}, E_{7}$ and $E_{8}$.\end{theorem}

\begin{remark} {\rm $r_2 G_2=3$ and $r_2 F_4=5$ were proved by Borel and  Serre in \cite{BS53}.}
\end{remark}

\section{Antipodal sets and real forms}\label{S8}

 Let $\psi$ be an involutive anti-holomorphic isometry of a Hermitian symmetric space  $M$ of compact type so that we have $\psi_{*}J=-J\psi_{*}$, where $J$ is the almost complex structure of $M$. Then the fixed point set $$F(\psi,M) = \{p \in M : \psi(p)=p\}$$ is called a {\it real form} of $M$ which is a connected totally geodesic Lagrangian submanifold $M$.
The classification of real forms of an irreducible Hermitian symmetric space of
compact type have been obtained by D. S. P Leung \cite{L79} and M. Takeuchi \cite{Ta84}. 

\vskip.05in
M. S. Tanaka and H. Tasaki proved in \cite{TT15} that a real form of a Hermitian
symmetric space $M$ of compact type is a product of real forms of irreducible factors
of $M$ and diagonal real forms determined from irreducible factors of $M$.

\vskip.06in
The following result implies that any two real forms in any Hermitian symmetric space of compact type have a non-empty intersection.

\begin{proposition} \cite{Cheng02,Ta10}. Let $M$ be a compact K\"ahler manifold whose holomorphic sectional curvatures are positive. If $L_{1}$ and $L_{2}$ are totally geodesic compact Lagrangian submanifolds of $M$, then $L_{1}\cap L_{2} \ne \emptyset$. \end{proposition}

The following theorem of M. Takeuchi  characterized real forms as symmetric $R$-spaces.

\begin{theorem} \cite{Ta84}  Every real form of a Hermitian symmetric space
of compact type is a symmetric $R$-space. Conversely, every symmetric $R$-space
is realized as a real form of a Hermitian symmetric space of compact type. The
correspondence is one-to-one.
\end{theorem}

 M. S. Tanaka and H. Tasaki  \cite{TT12} studied the intersection of two real forms in a Hermitian symmetric space of compact type. They proved  the following four results.
 
\begin{theorem}  \cite{TT12}  Let  $M$ be a Hermitian symmetric space of compact
type. If two real forms $L_1$ and $L_2$ of $M$ intersect transversally, then $L_1\cap  L_2$  is
an antipodal set of $L_1$ and $L_2$. \end{theorem}

\begin{theorem}  \cite{TT12}  Let $M$ be a Hermitian symmetric space of compact
type and let $L_1,L_2,L'_1,L'_2$ be real forms of $M$ such that $L_1,L'_1$ are congruent
and  $L_2,L'_2$ are congruent. If $L_1,L_2$ intersect transversally and if $L'_1,L'_2$
intersect transversally, then  $\#(L_1\cap  L_2) = \#(L'_1\cap  L'_2) $. 
 \end{theorem}
 
\begin{theorem}  \cite{TT12} Let $L_1,L_2$ be real forms of a Hermitian symmetric space  of compact type whose intersection is discrete. Then $L_1\cap L_2$ is an
antipodal set in $L_1$ and $L_2$. Moreover,   if $L_1$ and $L_2$ are congruent, then $L_1\cap L_2$ is a great antipodal set. Thus $\#(L_1\cap  L_2) = \#_2L_1 = \#_2L_2.$
\end{theorem} 

\begin{theorem} \cite{TT12} Let $M$ be an irreducible Hermitian symmetric
space of compact type and let $L_1,L_2$ be real forms of $M$ with $\#_{2}L_{1}\leq \#_{2}L_{2}$ and we assume that  $L_1\cap L_2$ is discrete. Then 
\vskip.06in
\begin{itemize}
\item[{\rm (a)}] If $M=G_{2m}(\mathbb C^{4m}) \, (m\geq 2)$, $L_{1}$ is congruent to $G_{m}(\mathbb H^{2m})$, $L_{2}$ is congruent to $U(2m)$, and
$$\#(L_1\cap  L_2) =2^{m} < \binom{2m}{m} =\#_2L_1<2^{2m} = \#_2L_2.$$
\item[{\rm (b)}] Otherwise, $\#(L_1\cap  L_2) = \#_2L_1 (\leq \#_2L_2).$
\end{itemize}\end{theorem} 

Y.-G. Oh defined in \cite{Oh91} the notion of global tightness of Lagrangian submanifolds in a Hermitian symmetric space; namely, a Lagrangian submanifold $L$ of a Hermitian symmetric space
$M$ is called {\it globally tight} if $L$ satisfies
$$\#(L\cap g\cdot L)=\dim H_{*}(L,\mathbb Z_{2})$$
for any isometry $g$ of $M$ such that $L$ intersects $g \cdot L$ transversally.

 \vskip.06in
H. Tasaki proved the following.

\begin{theorem} \cite{Ta10} In the complex hyperquadric, the intersection of two real
forms is an antipodal set whose cardinality attains the smaller 2-number of the two real forms. In particular,  every real form in the complex hyperquadric is a globally tight Lagrangian submanifold. \end{theorem}

\subsection{Fixed point sets and  holomorphic isometries}
In \cite{ITT15},   O. Ikawa, M. S. Tanaka and H. Tasaki discovered a necessary and sufficient condition for the fixed point set of a holomorphic isometry of a Hermitian symmetric space of compact type to be discrete. They also shown that the discrete fixed point set is an antipodal set. Further, they derived a necessary and sufficient condition that the intersection of two real forms in a Hermitian symmetric space of compact type is discrete. Moreover, they discussed some relations between the intersection of two real forms and the fixed point set of a certain holomorphic isometry by the use of the symmetric triads.

\begin{remark} {\rm For further results on real forms, we refer to \cite{Ik15,Ta13,Ta14}.}\end{remark}

  \section{Application to Lagrangian Floer homology}\label{S9}

Suppose that $(M,\omega)$ is a symplectic manifold, i.e., $M$ is a manifold equipped with a closed nondegenerate 2-form $\omega$. Let $L$ be a Lagrangian submanifold in $M$.
For a pair of closed Lagrangian submanifolds $(L_{0},L_{1})$ of $M$, one can define Lagrangian Floer homology $HF(L_{0},L_{1}:\mathbb Z_{2})$
with coefficient $\mathbb Z_{2}$ under some appropriate topological conditions.

\vskip.06in

 A. Floer \cite{Fl88} defined in 1988 the homology when $\pi_2(M,L_i)=0$, $i=0,1$. He
proved that it is isomorphic to the singular homology group $H_*(L_0,\mathbb Z_2)$ of $L_0$ in the case where $L_0$ is Hamiltonian isotopic to $L_1$. 
As a result,  Floer solved affirmatively the so called {\it Arnold conjecture} for Lagrangian intersections in that case (see \cite{Ar65,Fl88}). Symplectic Floer homology is invariant under Hamiltonian isotropy of the symplectomorphism.
Denote by $Hamilt(M,\omega)$ the set of all Hamiltonian diffeomorphisms of $M$.

\vskip.05in
In 1989,  A. Givental \cite{Gi89} proposed the following conjecture that generalized the results of Floer and himself.

\vskip.1in
\noindent {\bf Arnold-Givental Conjecture.} {\it Let $(M,\omega)$ be a symplectic manifold and
$\psi: M\to M$ be an anti-symplectic involution of $M$. Suppose that the fixed point set
$L=F(M,\psi)$ is compact and  nonempty. Then, for any $\phi\in Hamilt(M,\omega)$ such that the Lagrangian submanifold $L$ and its image $\phi(L)$ intersect transversally, the inequality
\begin{align}\label{10.1} \#(L\cap \phi(L))\geq \ b(L,\mathbb Z_2)\end{align}
holds, where $b(L,\mathbb Z_2)=\sum_{i\geq 0} b_i(L,\mathbb Z_2)$ is the total Betti number of $L$ with $\mathbb Z_{2}$ coefficient.}
\vskip.1in

 In  \cite{IST13}, H. Iriyeh, T. Sakai and H. Tasaki  computed   Lagrangian Floer homology $HF(L_0,L_1;{\mathbb Z}_2)$ for a pair of real forms 
$(L_0,L_1)$ in a monotone Hermitian symmetric space $M$ of compact type in the case where $L_0$ is not necessarily congruent to $L_1$. In particular, they established a generalization of the Arnold-Givental inequality \e{10.1} in the case where $M$ is irreducible. 
As an application, H. Iriyeh, T. Sakai and H. Tasaki established the following result.

\begin{theorem} \cite{IST13} Every totally geodesic Lagrangian sphere in the complex hyperquadric is globally volume minimizing under Hamiltonian deformations.
\end{theorem}

\section{Application to theories of designs and codes}\label{S10}
The theory of designs is the part of combinatorial mathematics that deals with the existence, construction and properties of systems of finite sets whose arrangements satisfy generalized concepts of balance and/or symmetry.

\vskip.1in
 \subsection{Codes and designs}\label{S10.1}

Codes and designs on association schemes are important research themes in combinatorics. In 1973, P. Delsarte  \cite{De73} gave linear programming bounds for cardinalities of codes and designs on commutative association schemes in terms of eigen-matrices. 

 \vskip.05in
After Delsarte's work, the theory of designs on spheres was introduced  in 1977 by P. Delsarte, J. M. Goethals and J. J. Seide in \cite{DGS77} as an analogy of Delsarte theory. The main tool in their works is the addition formula for polynomials; polynomials associated with metric or cometric association schemes, or the Gegenbauer polynomials with spheres. 
As a result, Delsarte's bounds  were established in terms of spherical Fourier transforms. For a survey on the studies of codes and designs on spheres, we refer to \cite{BB09}.

\vskip.05in
Compact symmetric spaces of rank one are natural and significant examples of the Delsarte spaces or the polynomial spaces for continuous metric spaces. The theory of designs on rank one compact symmetric spaces was also investigated by S. G. Hoggar \cite{Ho82} in details.
E. Bannai and S. G. Hoggar  also studied on rank one compact symmetric spaces in \cite{BH85}. 

\vskip.05in
For other compact symmetric spaces, codes, designs and Delsarte's bounds have been studied by many researchers. For examples, studied by C. Bachoc, R. Coulangeon and G. Nebe \cite{BCN02} and  C. Bachoc, E. Bannai and R. Coulangeon \cite{BBC04} on real Grassmannian manifolds; by A. Roy \cite{R10} on complex Grassmannian manifolds; and by A. Roy and A. J. Scott \cite{RS09} on unitary groups.

\vskip.05in
In \cite{KO20.1}, H. Kurihara and T. Okuda provided a definition of codes and designs  on general compact symmetric spaces. They also established in \cite{KO20.1} a general formulation of Delsare's bounds on compact symmetric spaces.

\vskip.1in
 \subsection{Great antipodal sets and designs on complex Grassmannian manifolds}\label{S10.2}

For compact symmetric spaces of higher rank, H. Kurihara and T. Okuda \cite{KO20}  obtained a characterization of maximal antipodal sets of complex Grassmannians in term of certain designs (more precisely, ${\mathcal E}\cup {\mathcal F}$-designs) with the smallest cardinalities. 
In particular, Kurihara and Okuda's  main result in \cite{KO20}  implies the following.

\begin{theorem} \cite{KO20} A great antipodal set of a complex  Grassmannian manifold is an ${\mathcal E}$-design with the smallest cardinality.
\end{theorem}

 \subsection{Cubature formulas for great antipodal sets on complex Grassmannian manifolds}\label{S10.3}

In \cite{O17}, H. Kurihara  and T. Okuda established a formulation of Delsarte theory for finite subsets of compact symmetric spaces. As its application, they proved that great
antipodal subsets of complex Grassmannian manifolds give rise to  cubature
formulas for certain functional spaces.

\vskip.1in
 \subsection{Great antipodal sets and designs on unitary groups}\label{S10.4}
 
Put  $[n]=\{1, 2, . . . , n\}$ and let $2^{[n]}$ denote the power set of $[n]$. Let $Q$
be the $n$-ary Cartesian product of the two elements set $\{1,-1\}$. Then the {\it Hamming cube
graph} $Q_{n}$ of degree $n$ is the graph with the vertex set $Q$ and two vertices are adjacent
whenever they differ in precisely one coordinate.

\vskip.05in
In \cite{Ku21} H. Kurihara investigated a relation between great antipodal sets on unitary group $U(n)$ and design theory on $U(n)$. In \cite{Ku21}, he also established a beautiful relationship between a great antipodal set on $U(n)$ and a Hamming cube graph $Q_{n}$.

\subsection{Application to coding theory}\label{S10.5}
Coding theory studies properties of codes and their respective fitness for specific applications. 
The main purpose of codes is to be able to recover the original content of a transmitted message by correcting errors that have entered the message during transmission. This capability is useful in maintaining the integrity of computer networks, communication systems, compact disk recording, etc.

\vskip.05in
A $p$-group $H$ is called {\it extra special} if its center ${\mathbb Z}$ is cyclic of order $p$, and the quotient $H/{\mathbb Z}$ is a non-trivial elementary abelian $p$-group.
In 1989, J. A. Wood \cite{Wo89} investigated the equivalence between the diagonal extra-special 2-group  of spinor $Spin(n)$ and the self-orthogonal linear binary codes of algebraic coding theory.  In Wood's article \cite{Wo89}, Theorem \ref{T:7.8} and Theorem \ref{T:7.9} were mentioned and used.

\section{$k$-symmetric spaces, $\Gamma$-symmetric spaces, $k$-number and flag manifolds}\label{S11}

\subsection{$k$-symmetric spaces and $\Gamma$-symmetric spaces} \label{S11.1}

Since the 1960s, generalizations of symmetric spaces have been proposed in various directions. In 1967, A. J. Ledger \cite{Le67} initiated the study of  $s$-manifolds. These are Riemannian manifolds $M$ which admit at each point $x\in M$ a symmetry $s_{x}$ with $x$ as an isolated fixed point. A $k$-symmetric structure is called {\it regular} \cite{Ko80} if it satisfies 
 \begin{align}\label{12.1} \theta_x\circ \theta_y=\theta_z\circ \theta_z,\;\; z=\theta_x(y).\end{align} 
If $s_{x}$ is
of finite order $k$, a regular $s$-manifold is called a $k$-symmetric space (see \cite{Ko77}).

\vskip.05in
As a further generalization of Riemannian symmetric spaces,  P. Lutz \cite{Lu81} introduced in 1981 $\Gamma$-symmetric space, where  $\Gamma$ is a finite abelian group. These are manifolds $M$ which admit the following structure: To each point $x\in M$ one assigns in a suitable way a group $\Gamma_{x}$ isomorphic to $\Gamma$ which acts effectively on $M$ with $x$ as an isolated fixed point. If $\Gamma$ is isomorphic to $\mathbb Z_{2}$, then a $\Gamma$-symmetric space is just a Riemannian symmetric space.
 
\vskip.05in
Every complex flag manifold can be regarded as an $R$-space.
 Let $U$ be a compact connected semisimple centerless Lie group and let $\frak u$ be the Lie algebra of $U$. Then the {\it complex flag manifold} of $U$ is the orbit of the adjoint action of $U$ on $\frak u$. 
Take $M= Ad(U)Y$ for $Y \ne 0$ in $U$ and let ${\mathfrak g}= {\mathfrak u}_{\mathbb C}= {\mathfrak u}+i {\mathfrak u}$. Then there exists a Cartan decomposition of the realization ${\mathfrak g}_{\mathbb R}$  of ${\mathfrak g}$ and one may consider $M$ as the orbit of $iY$  in $i{\mathfrak u}$ by the adjoint action of $U$.

\vskip.1in

\subsection{Maximal antipodal sets and  $\Gamma$-symmetric $R$-spaces}\label{S11.2}

In 2020, P. Quast and T. Sakai \cite{QS20} extended the definition of antipodal sets of compact symmetric spaces naturally extends to $\Gamma$-symmetric spaces as follows (see also \cite{QS22}).

\begin{definition} \cite{QS22} {\rm Let $\Gamma$ be a finite abelian group, and let $\mu=\{\mu^{\gamma}\}_{\gamma\in  \Gamma}$ be a
$\Gamma$-symmetric structure on a manifold $M$. A subset $A$ of the $\Gamma$-symmetric space $(M,\mu)$ is called antipodal if $\gamma_{x}(y)=y$ for all $x,y\in A$ and for all $\gamma\in \Gamma$.
An antipodal set $A$ of $(M,\mu)$ is called {\it maximal} if $A$ is not a proper subset of another
antipodal set of $M$. The supremum of the cardinalities of antipodal sets of $(M,\mu)$
is called the antipodal number denoted by $\#_{\Gamma}M$. An antipodal set $A$ of $(M,\mu)$ is called {\it great} if the cardinality of $A$ is equal to $\#_{\Gamma}M$.}
\end{definition}

In \cite{QS20}, P. Quast and T. Sakai defined the induced natural $\Gamma$-symmetric structure on $R$-spaces. Further, they determined the maximal antipodal sets of $R$-spaces with respect to the induced natural $\Gamma$-symmetric structures. In particular, they shown that
 any two maximal antipodal sets of a $R$-space with respect to an induced natural $\Gamma$-symmetric structures are conjugate.

\vskip.1in

\subsection{A link between real flag manifolds and complex flag manifolds}\label{S11.3} 

 In 1997, C. U. S\'anchez proved the following.

\begin{proposition}\label{P:12.1} \cite{Sa97} If $M$ is a real flag manifold, then there exists a complex flag manifold $M_{\mathbb C}$ such that $M$ is isometrically imbedded in $M_{\mathbb C}$. If $M$ is a symmetric $R$-space, then $M_{\mathbb C}$ is a hermitian symmetric space and the isometric imbedding is totally geodesic. If $M$ is already a complex flag manifold, then $M_{\mathbb C} = M$.
\end{proposition}

\begin{remark} {\rm  H. Iriyeh, T. Sakai and H. Tasaki \cite{IST14,IST19} proved that the intersection of real flag manifolds in the complex flag manifold consisting of sequences of complex subspaces in a complex vector space is an antipodal set, which is a generalization of that in a Hermitian symmetric space of compact type.}
\end{remark}

\vskip.1in
 \subsection{$k$-number, index number and complex flag manifold}\label{S11.4} 
 
 For a complex flag manifold $M_{\mathbb C}$, there exists a positive integer $k_{0}=k_{o}(M_{\mathbb C})\geq 2$ such that, for each integer $k\geq k_{0}$, there exists a {\it $k$-symmetric structure} \cite{Ko80} on $M_{\mathbb C}$, i.e., for each point $x\in M_{\mathbb C}$ there exists an isometry $\theta_x$ such that $\theta_x^k=id$ with $x$ as an isolated fixed point.

\vskip.06in
Analogous to Proposition \ref{P:2.9} for 2-number of  compact symmetric spaces,  C. U. S\'anchez defined {\it $k$-number}, denoted by $\#_{k}(M_{\mathbb C})$, of a complex flag manifold $M_{\mathbb C}$ as the maximal possible cardinality of the $k$-sets $A_{k}\subset M_{\mathbb C}$ which satisfies the property that for each point $x\in A_{k}$ the corresponding $k$-symmetry at $x$ fixes every point in $A_{k}$. 
 
 \vskip.1in
As an extension of Theorem \ref{T:4.16} of Takeuchi, C. U. S\'anchez proved the following. 

\begin{theorem} \cite{Sa93}  For each complex flag manifold $M_{\mathbb C}$, we have
$ \#_{k}(M_{\mathbb C})=\dim H^{*}(M_{\mathbb C},\mathbb Z_{2})$,

\end{theorem}

Applying Proposition \ref{P:12.1}, S\'anchez  \cite{Sa97} defined in 1997  the {\it index number} of a real flag manifold $M$, denoted by $\#_{I}M$, as the maximal possible cardinality of the $p$-sets $A_{p}M$ with $p$ a prime number, in terms of fixed points of symmetries of the complex flag manifolds restricted to the real one.   

 \vskip.1in
  
For index number of a real flag manifold, C. U. S\'anchez obtained the following result in 1997.
 
 \begin{theorem} \cite{Sa97} Let $M$ be a real flag manifold. Then  $\#_{I}M=b(M,\mathbb Z_2).$ 
 \end{theorem}
 
 \vskip.1in
  \subsection{$k$-number and generalized flag manifolds}\label{S11.5}
  
 Let $G$ be a compact connected semisimple Lie group. Then the homogeneous spaces one obtains as orbits of $G$ under the adjoint representation on the Lie algebra of $G$ are also called {\it generalized flag manifolds}. 
It is known that every generalized flag manifold admits $k$-symmetric structure.

In a similar way, S\'anchez \cite{Sa97} also proved the following.

\begin{theorem}  If $M$ is a  generalized flag manifold, then  $\#_{k}(M)=\chi(M)$
for any  $k$-symmetric structure on $M$.
\end{theorem}

 \vskip.1in
 \subsection{$k$-number and generalized flag manifolds}\label{S11.6}
 
 Let $G$ be a compact connected semisimple Lie group. Then the homogeneous spaces one obtains as orbits of $G$ under the adjoint representation on the Lie algebra of $G$ are also called {\it generalized flag manifolds}. 
It is known that every generalized flag manifold admits $k$-symmetric structure.

  \vskip.05in
In a similar way, C. U. S\'anchez also proved the following.

\begin{theorem} \label{T:11.3} If $M$ is a  generalized flag manifold, then any of its $k$-symmetric structure on $M$ satisfies $\#_{k}(M)=\chi(M).$
\end{theorem}

 \vskip.1in
 \subsection{$k$-number and $k$-symmetric submanifolds}\label{S11.7}
 
Let $M\subset {\mathbb E}^m$ be submanifold of ${\mathbb E}^m$.  If $M$ satisfies   \begin{itemize}
\vskip.02in
  \item[{\rm (a)}]   For each $x\in M$, there is an isometry $\sigma_x: {\mathbb E}^m\to {\mathbb E}^m$ such that 
  $\sigma_x^k =id_M$, $\sigma_x(x)=x$, and $\sigma_x |_{T^\perp_x M}$ = identity on $T^\perp_x M$; 
   
 \vskip.02in
  \item[{\rm (b)}]  $\sigma_x(M)\subset M$;  and 
   
  \vskip.02in \item[{\rm (c)}]  Let $\theta_x=\sigma_x |_M$. The collection $\{\theta_x, \,x\in M\}$ defines on $M$ a 
    Riemannian regular $s$-structure of order $k$,
      \end{itemize}
       then $M$ is called  an {\it  extrinsic $k$-symmetric submanifold} (see \cite{Fe80}).
 
C. U. S\'anchez  \cite{Sa93} proved the following result.
 
\begin{theorem} If $M\subset \mathbb E^m$ is an  extrinsic $k$-symmetric submanifold, then 
 $ \#_{k}(M)=b(M,\mathbb Z_p)$
 for any  prime number $p\geq 2$ which divides $k$.
 \end{theorem}

  \vskip.1in
 \subsection{Morse functions and great antipodal sets on $G_{2}/SO(4)$}\label{S11.8} 

  In \cite{Sa22.3}, Y. Sasaki constructed $\mathbb Z_{2}$-perfect Morse functions of $GI=G_{2}/SO(4)$ whose set of all critical points is a great antipodal set of  $GI$. Consequently, he provided a reason why the 2-number $\#_{2}(GI)$ matches the Betti-number of the $Z_{2}$-coefficient homology group of  $GI$.

 \section{2-number, index number and CW complex structure} \label{S12}

The following conjecture was posed first time in author's 1987 report \cite{C87}.
 
 \vskip.05in
 \noindent {\bf Conjecture 1. } {\it  For any a compact symmetric space $M$, $\#_{2}M$ is equal to the smallest number of cells that are needed for a CW complex structure on $M$.}
  \vskip.1in

Related to this conjecture, J. Berndt, S. Console and A. Fino proved the following.
 
 \begin{theorem} \cite{BCF01} The index number $\#_{I}M$ is equal to the smallest number of cells  that are needed for a CW complex structure  for each real flag manifold $M$.
 \end{theorem}
 
 In the proof of this theorem, the authors have applied the convexity theorems of M. F.  Atiyah \cite{At82}, V. Guillemin and S. Sternberg's result in \cite{GS82} for symplectic manifolds with a hamiltonian torus action as well as a generalization of J. J. Duistermaat's result in \cite{Du83} for fixed point set of antisymplectic involutions.

The next result was also proved by Berndt, Console and Fino in \cite{BCF01}.

 \begin{theorem} The index number $\#_{I}M$ of a real flag manifold $M$ satisfies
 $\#_{I}M=\chi(M)\, ({\rm mod\, 2}).$
  \end{theorem}

 In \cite{CN88}, the author and T. Nagano made the following:
 
\vskip.07in
 \noindent {\bf Conjecture 2. } {\it  $\#_{2}M=\chi(M)\, ({\rm mod  \,2})$ holds for every  irreducible compact symmetric space $M$.}

 \vskip.07in
It was known that the total Betti numbers of a simply-connected compact symmetric space $M$ satisfies $b(M;{\bf R}) \leq \#_2M$ (see \cite[page 54]{C87}).
 Professor T. Nagano asked  the following open problem.

  \vskip.1in
\noindent {\bf Problem.} $b(M;{\bf R})< \#_2M$  $ \Longrightarrow$   $M$ has 2-torsion?
 
  \vskip.1in
  As far as I know, this problem remains open till now.


\begin{thebibliography}{106}

\bibitem{Ad87} J. F. Adams, {\it 2-tori in $E_{8}$}, {Math. Ann.} {\bf 278} (1987), 29--39.

\bibitem{Ar65} V. I. Arnold, {\it Sur une propri\'et\'e topologique des applications globalement canoniques de la m\'ecanique classique}, C. R. Acad. Sci. Paris {\bf 261} (1965), 3719--3722.


\bibitem{At61} M. F. Atiyah, {\it Characters and cohomology of finite groups}, {Inst. Hautes \'Etudes Sci. Publ. Math}. No. {\bf 9} (1961), 23--64. 

\bibitem{At82} M. F. Atiyah, {\it  Convexity and commuting hamiltonians}, Bull. London Math. Soc. {\bf 14} (1982), 1--15.

\bibitem{ABS64} M. F. Atiyah, R. Bott and A. Shapiro, {\it Clifford modules}, {Topology} {\bf 3},  (1964), Suppl. 1, 3--38.

\bibitem{AS68} M. F. Atiyah and I. M. Singer, {\it The index of elliptic operators III}, Ann. of Math. {\bf 87} (1968), 546--604.

\bibitem{BBC04}  C. Bachoc, E. Bannai and R. Coulangeon, {\it Codes and designs in Grassmannian spaces}, Discrete Math. {\bf 277} (2004), 15--28.


\bibitem{BCN02} C. Bachoc, R. Coulangeon and G. Nebe, {\it Designs in Grassmannian spaces and lattices}, J. Algebraic Combin. {\bf 16} (2002), 5--19.

\bibitem{BB09} E. Banna and E. Bannai, {\it A survey on spherical designs and algebraic combinatorics on spheres}, European J. Combin. {\bf 30} (2009), 1392--1425.

\bibitem{BH85} E. Bannai and S. G. Hoggar, {\it On tight $t$-designs in compact symmetric spaces of rank one}, Proc. Japan Acad. Ser. A Math. Sci. {\bf  61}, (1985), 78--82, (1985).

\bibitem{BCF01}  J. Berndt, S. Console and A. Fino, {\it On index number and topology of flag manifolds}, Differential Geom. Appl. {\bf 15(1)} (2001),  81--90. 

\bibitem{Bey18} J. Beyrer, {\it  A complete description of the antipodal set of most symmetric spaces of compact type}, Osaka J. Math. {\bf 55(3)}  (2018), 567--586.

\bibitem{BS53} A. Borel and J.-P. Serre, {\it Sur certains sousgroupes des groupes de Lie compacts}, {Comm. Math. Helv.} {\bf 27} (1953), 128--139.

\bibitem{B33} K. Borsuk, {\it  Drei Satze uber die $n$-dimensionale euklidische Sph\"are}, Fund. Math. {\bf 20} (1933), 177--190.

\bibitem{Bo59} R. Bott, {\it  Stable homotopy of the classical groups}, Ann. of Math. {\bf 70} (1959), 313--337.

\bibitem{Bo70} R. Bott, {\it The periodicity theorem for the classical groups and some of its applications},  Adv. Math. {\bf 4} (1970), 353--411.

\bibitem{Bu92} J. M. Burns, {\it Homotopy of compact symmetric spaces}, Glasgow Math. J. {\bf 34(2)} (1992),  221--228.

\bibitem{C73} B.-Y. Chen, {\it  Geometry of submanifolds}, Marcel Dekker, New York, NY, 1973.

\bibitem{C87}  B.-Y. Chen,  {\it A new approach to compact symmetric spaces and applications},  (A report on joint work with Professor T. Nagano), Katholieke Universiteit Leuven, Belgium, 1987.

\bibitem{C87.2}  B.-Y. Chen,  {\it Symmetries of compact symmetric spaces}, Geometry and Topology of Submanifolds (Marseille, 1987), pp. 38--48, World Scientific, Teaneck, NJ, 1989.

\bibitem{C13}  B.-Y. Chen,  {\it The 2-ranks of connected compact Lie groups}, Taiwanese J. Math. {\bf 17(3)} (2013),  815--831.

\bibitem{C17}   B.-Y. Chen,  {\it Borsuk-Ulam theorem and maximal antipodal sets of compact symmetric spaces}, Int. Electron. J. Geom. {\bf 10(2)} (2017),  11--19.

\bibitem{C18}   B.-Y. Chen,  {\it Two-numbers and their applications--a survey},
Bull. Belg. Math. Soc. Simon Stevin {\bf 25(4)} (2018),  565--596.

\bibitem{CLN80} B.-Y. Chen, P.-F. Leung and T. Nagano, {\it Totally geodesic submanifolds of symmetric spaces  III}, preprint, 1980. DOI: 10.48550/arXiv.1307.7325

\bibitem{CN77}  B.-Y. Chen and T. Nagano,  {\it Totally geodesic submanifolds of symmetric spaces I}, Duke Math. J. {\bf 44} (1977), 745--755.

\bibitem{CN78} B.-Y. Chen and T. Nagano, {\it Totally geodesic submanifolds of symmetric spaces II}, Duke Math. J. {\bf 45} (1978),  405--425.

\bibitem{CN82} B.-Y. Chen and T. Nagano,  {\it Un invariant g\'eom\'etrique riemannien}, C. R. Acad. Sci. Paris S\'er. I Math. {\bf 295(5)} (1982),  389--391.

\bibitem{CN88}  B.-Y. Chen and T. Nagano, {\it A Riemannian geometric invariant and its applications to a problem of Borel and Serre}, Trans. Amer. Math. Soc. {\bf 308} (1988),  273--297. 

\bibitem{Cheng02}  X. Cheng, {\it The totally geodesic coisotropic submanifolds in K\"ahler manifolds},  Geom. Dedicata {\bf 90} (2002), 115--125.

\bibitem{De73} P. Delsarte,  {\it An algebraic approach to the association schemes of coding theory}, Philips Res. Rep. Suppl., {\bf 10} (1973).

\bibitem{DGS77} P. Delsarte, J. M. Goethals and J. J. Seidel, {\it Spherical codes and designs},
Geom. Dedicata, {\bf 6(3)} (1977),  363--388.

\bibitem{Du83} J. J. Duistermaat, {\it Convexity and tightness for restrictions of Hamiltonian functions to fixed point sets of an antisymplectic involution}, Trans. Amer. Math. Soc. {\bf 275} (1983) 417--429.

\bibitem{EM62} S. Eilenberg and J. C. Moore, {\it Limits and spectral sequences}, {Topology} {\bf 1} (1962), 1--23.

\bibitem{EQT15} J.-H. Eschenburg, P. Quast and M. S. Tanaka, {\it Maximal tori of extrinsic symmetric spaces and meridians},  Osaka J. Math. {\bf 52(2)} (2015), 299--305.

\bibitem{Fa93} A. Fauntleroy, {\it Projective ranks of Hermitian symmetric spaces}, Math. Intelligencer {\bf 15} (1993),  27--32.

\bibitem{Fe80} D. Ferus, {\it Symmetric submanifolds of Euclidean space}, Math. Ann. {\bf 247(1)} (1980), 81--93.

\bibitem{Fl88} A. Floer, {\it Morse theory for Lagrangian intersections}, J. Differential Geom. {\bf 28} (1988), 513--547.

\bibitem{Gi89}  A. Givental, {\it Periodic maps in symplectic topology}, Funct. Anal. Appl. {\bf  23} (1989), 37--52.

\bibitem{GS82} V. Guillemin and S. Sternberg, {\it Convexity properties of the moment mapping}, Invent. Math. {\bf 67} (1982), 491--513.

\bibitem{H} S. Helgason, {\it Differential geometry, Lie groups and symmetric spaces}, Academic Press, New York, 1978.

\bibitem{Ho82} S. G. Hoggar, {\it t-designs in projective spaces}, European J. Combin. {\bf 3(3)} (1982),  233--254.

\bibitem{HS} H. Hopf and H. Samelson, {\it Ein Satz uber Wirkungsraume geschlossener Liescher Gruppen}, {Comm. Math. Helv.} {\bf13} (1941), 240--251.

\bibitem{Ik15} O. Ikawa, H. Iriyeh, T. Okuda, T. Sakai and H. Tasaki, {\it
Antipodal structure of the intersection of real flag manifolds in a complex flag manifold II},
Proceedings of the 19th International Workshop on Hermitian-Grassmannian Submanifolds and 10th RIRCM-OCAMI Joint Differential Geometry Workshop, pp. 147--157, Natl. Inst. Math. Sci. (NIMS), Taej\u{o}n, 2015.

\bibitem{ITT15} O. Ikawa, M. S. Tanaka and H. Tasaki, {\it The fixed point set of a holomorphic isometry, the intersection of two real forms in a Hermitian symmetric space of compact type and symmetric triads}, Internat. J. Math. {\bf 26(6)} (2015),  1541005, 32 pp.

\bibitem{IST13} H. Iriyeh, T. Sakai and H. Tasaki, {\it Lagrangian Floer homology of a pair of real forms in Hermitian symmetric spaces of compact type},  J. Math. Soc. Japan {\bf 65} (2013),  1135--1151. 

\bibitem{IST14} H. Iriyeh, T. Sakai and H. Tasaki, {\it Lagrangian intersection theory
and Hamiltonian volume minimizing problem}, Real and Complex submanifolds, pp. 391--399.
Springer Proc. Math. Stat., {\bf 106}, Springer, Tokyo, 2014.

\bibitem{IST19} H. Iriyeh, T. Sakai and H. Tasaki, {\it On the structure of the intersection of real flag manifolds in a complex flag manifold}, Differential geometry and Tanaka theory--differential system and hypersurface theory, pp. 87--98. Adv. Stud. Pure Math., {\bf 82} Mathematical Society of Japan, Tokyo, 2019.

\bibitem{KNa64} S. Kobayashi and T. Nagano, {\it On filtered Lie algebras and geometric structures I}, J. Math. Mech. {\bf 13} (1964), 875--907.

\bibitem{Ko77} A. Kono, {\it On the 2-rank of compact connected Lie groups}, {J. Math. Kyoto Univ.} {\bf 17} (1977), 1--18.

\bibitem{Ko80} O.  Kowalski, {\it Generalized symmetric spaces}, Lecture Notes in Math. volume {\bf 805}, Springer, Berlin-Heidelberg-New York, 1980.

\bibitem{Ku21} H. Kurihara, {\it Antipodal sets and designs on unitary groups},
Graphs Combin. {\bf 37(5)} (2021),  1559--1583.

\bibitem{KO20} H. Kurihara and T. Okuda, {\it Great antipodal sets of complex Grassmannian manifolds as designs with the smallest cardinalities}, J. Algebra {\bf 559} (2020), 432--466.

\bibitem{KO20.1} H. Kurihara and T. Okuda, {\it Combinatorics and Fourier analysis on compact symmetric spaces}, Quandles and Symmetric Spaces, pp.61--73, OCAMI Report, {\bf 4} (2021).
doi: 10.24544/ocu.20210605-004

\bibitem{Le67} A. J. Ledger, {\it Espaces de Riemann sym\'etriques g\'en\'eralis\'es},  C. R. Acad. Sci. Paris S\'er. A-B {\bf 264} (1967), A947--A948.

\bibitem{L79} D. P. S. Leung, {\it Reflective submanifolds. IV, Classification of real forms of Hermitian symmetric spaces}, J. Differential Geometry {\bf 14(2)} (1979), 179--185.

\bibitem{LD14} X. Liu and S. Deng, {\it The antipodal sets of compact symmetric spaces}, Balkan J. Geom. Appl. {\bf 19(1)}, (2014), 73--79.

\bibitem{Lo85} O. Loos, {\it Charakterisierung symmetrischer R-R\"aume durch ihre Einheitsgitter}, Math. Z. {\bf 189(2)} (1985), 211--226. 

\bibitem{Lu81} R. Lutz, {\it Sur la g\'eom\'etrie des espaces $\Gamma$-sym\'etriques},
C. R. Acad. Sci. Paris S\'er. I Math. {\bf 293(1)} (1981),  55--58

\bibitem{LF51}  L. A. Lyusternik and A. I. Fet, {\it Variational problems on closed manifolds},  Doklady Akad. Nauk SSSR (N.S.) {\bf 81} (1951). 17--18. 

\bibitem{Ma65} P.  May, {\it The cohomology of restricted Lie algebras and of  Hopf algebras}, {Bull. Amer. Math. Soc.} {\bf 71} (1965), 372--377.

\bibitem{MM65} J. Milnor and J. Moore, {\it On the structure of Hopf algebras}, {Ann. of Math.} {\bf 81} (1965), 211--294.

\bibitem{Na65} T. Nagano, {\it  Transformation groups on compact symmetric spaces}, Trans. Amer. Math. Soc. {\bf 118} (1965), 428--453.

\bibitem{Na88} T. Nagano,{\it  The involutions of compact symmetric spaces}, {Tokyo J. Math.} {\bf 11} (1988),  57--79. 

\bibitem{Na92} T. Nagano, {\it The involutions of compact symmetric spaces II}, {Tokyo J. Math.} {\bf 15} (1992),  39--82. 

\bibitem{NM88} T. Nagano and M. Sumi, {\it The structure of the symmetric space with applications},  Geometry of manifolds (Matsumoto, 1988), pp. 111--128,  Perspect. Math., 8, Academic Press, Boston, MA, 1989.

\bibitem{Oh91} Y.-G. Oh, Tight Lagrangian submanifolds in $CP^{n}$, Math. Z. {\bf 207} (1991), 409--416.

\bibitem{O17} T. Okuda, {\it Cubature formulas for great antipodal sets on complex Grassmann manifolds}, Combinatorics of Lie Type, RIMS Kokyuroku {\bf 2039} (2017),  79--89. https://core.ac.uk/download/pdf
/188809952.pdf

\bibitem{Qu13} P. Quast, {\it Centrioles in symmetric spaces}, Nagoya Math. J. {\bf 211} (2013), 51--77.

\bibitem{QT12} P. Quast and M. S. Tanaka, {\it Convexity of reflective submanifolds in symmetric R-spaces}, Tohoku Math. J.  {\bf 64(4)} (2012),  607--616. 

\bibitem{QS20} P. Quast and T. Sakai, {\it  Natural $\Gamma$-symmetric structures on $R$-spaces}, J. Math. Pures Appl. (9) {\bf 141} (2020), 371--383.

\bibitem{QS22} P. Quast and T. Sakai, {\it A survey on natural $\Gamma$-symmetric structures on R-spaces}, Differential geometry and global analysis--in honor of Tadashi Nagano, Contemp. Math. {\bf 777} (2022), 185--197. 

\bibitem{Qu71.1} D. Quillen, {\it The spectrum of an equivariant cohomology ring}, {Ann. of Math.} {\bf 94} (1971), 549--409.

\bibitem{R10} A. Roy, {\it Bounds for codes and designs in complex subspaces}, J. Algebraic Combin. {\bf 31} (2010), 1--32.

\bibitem{RS09} A. Roy and A. J. Scott, {\it Unitary designs and codes}, Des. Codes Cryptogr. {\bf 53} (2009), 13--31.

\bibitem{Sa93} C. U. S\'anchez, {\it  The invariant of Chen-Nagano on flag manifolds}, Proc. Amer. Math. Soc. {\bf 118(4)} (1993),  1237--1242. 

\bibitem{Sa97} C. U. S\'anchez, {\it The index number of an R-space: an extension of a result of M. Takeuchi's},  Proc. Amer. Math. Soc. {\bf 125(3)} (1997),  893--900. 

\bibitem{SG02} C. U. S\'anchez and A. Giunta,  {\it The projective rank of a Hermitian symmetric space: a geometric approach and consequences}, Math. Ann. {\bf 323(1)} (2002), 55--79.

\bibitem{Sa21} Y. Sasaki, {\it  A study on maximal antipodal sets of compact symmetric spaces}, 
PhD Thesis, University of Tsukuba, Tokyo, 2021.

\bibitem{Sa22} Y. Sasaki, {\it  Homogeneity of maximal antipodal sets},  Osaka J. Math. {\bf 59(1)} (2022),  115--144.

\bibitem{Sa22.1} Y. Sasaki, {\it  Maximal antipodal sets of  $F_{4}$ and  $FI$ },
J. Lie Theory {\bf 32(1)} (2022),  281--300.

\bibitem{Sa22.2} Y. Sasaki, {\it Maximal antipodal sets of  $E_{6}$  and some compact symmetric spaces}, Differential Geom. Appl. {\bf 85} (2022), Paper No. 101934, 29 pp.

\bibitem{Sa22.3} Y. Sasaki, {\it  Morse Functions of $G_{2}/SO(4)$,}
Tokyo J. Math. {\bf 45(1)} (2022),  201--214.

\bibitem{S85} H. Steinlein, {\it  Borsuk's antipodal theorem and its generalizations and applications: a survey},  S\'em. Math. Sup\'er. Montr\'eal, S\'em.  {\bf 95} (1985), 166--235.

\bibitem{Sw71} R. G. Swan, {\it Groups with no odd dimensional cohomology}, {J. Algebra} {\bf 17} (1971) 401--403. 

\bibitem{Ta65} M. Takeuchi, {\it Cell decompositions and Morse equalities on certain symmetric spaces}, {J. Fac. Sci. Univ. Tokyo} {\bf  12} (1965), 81--192.

\bibitem{Ta84} M. Takeuchi, {\it Stability of certain minimal submanifolds of compact Hermitian symmetric spaces}, Tohoku Math. J. {\bf 36} (1984), 293--314.

\bibitem{Ta89} M. Takeuchi,  {\it Two-number of symmetric $R$-spaces}, {Nagoya Math. J.} {\bf 115} (1989), 43--46. 

\bibitem{Ta11} M. S. Tanaka, {\it  Antipodal sets of compact Riemannian symmetric spaces and their applications}, Proceedings of the 15th International Workshop on Differential Geometry and the 4th KNUGRG-OCAMI Differential Geometry Workshop, {\bf 15}, pp. 11--18.
National Institute for Mathematical Sciences (NIMS), Taej\u{o}n, 2011

\bibitem{Ta13} M. S. Tanaka, {\it Fixed point sets of isometries and the intersection of real forms in a Hermitian symmetric space of compact type}, Proceedings of The Seveneenth International Workshop on Diff. Geom. {\bf 17} (2013), pp. 1--9.

\bibitem{Ta14} M. S. Tanaka, {\it  Antipodal sets of compact symmetric spaces and the
intersection of totally geodesic submanifolds}, Differential Geometry of Submanifolds and its Related Topics, pp. 205--219 (2014), World Scientific, Hackensack, NJ. 

\bibitem{Ta15} M. S. Tanaka,  {\it Geometry of symmetric R-spaces}, Geometry and analysis on manifolds, 471--481, Progr. Math., 308, Birkhauser/Springer, 2015.  

\bibitem{TT12} M. S. Tanaka and H. Tasaki,  {\it The intersection of two real forms in Hermitian symmetric spaces of compact type}, J. Math. Soc. Japan {\bf 64(4)} (2012), 1297--1332. 

\bibitem{TT13} M. S. Tanaka and H. Tasaki, {\it Antipodal sets of symmetric R-spaces}, Osaka J. Math. {\bf 50(1)} (2013), 161--169.

\bibitem{TT15} M. S. Tanaka and H. Tasaki,  {\it The intersection of two real forms in Hermitian symmetric spaces of compact type II}, J. Math. Soc. Japan {\bf 67(1)} (2015),  275--291. 

\bibitem{TT15.1} M. S. Tanaka and H. Tasaki, {\it  Correction to: ``The intersection of two real forms in Hermitian symmetric spaces of compact type''}  J. Math. Soc. Japan {\bf 67(3)} (2015),  1161--1168.

\bibitem{TT17} M. S. Tanaka and H. Tasaki, {\it Maximal antipodal subgroups of some compact classical Lie groups},  J. Lie Theory {\bf 27} (2017), 801--829.

\bibitem{TT20} M. S. Tanaka and H. Tasaki, {\it  A Maximal antipodal sets of compact classical symmetric spaces and their cardinalities I}, Differential Geom. Appl. {\bf 73} (2020), 101682, 32 pp.

\bibitem{TTY19} M. S. Tanaka, H. Tasaki and O. Yasukura, {\it 
Maximal antipodal sets of $G_{2}$ and $G_{2}/SO(4)$ and related geometry}, Proceedings of the 22nd International Workshop on Differential Geometry of Submanifolds in Symmetric Spaces \& the 17th RIRCM-OCAMI Joint Differential Geometry Workshop, pp. 153--159 (2019). 

\bibitem{TTY22} M. S. Tanaka, H. Tasaki and O. Yasukura, {\it  Maximal antipodal sets related to  $G_{2}$,} Proc. Amer. Math. Soc. {\bf 150(10)} (2022),  4533--4542. 

\bibitem{Ta10} H. Tasaki, {\it The intersection of two real forms in the complex hyperquadric}, Tohoku Math. J. {\bf 62(3)} (2010),  375--382. 

\bibitem{Tas13} H. Tasaki,  {\it Antipodal sets in oriented real Grassmann manifolds}, Internat. J. Math. {\bf 24(8)} (2013),  1350061, 28 pp.

\bibitem{Tas14} H. Tasaki,  {\it Sequences of maximal antipodal sets of oriented real Grassmann manifolds}, Real and complex submanifolds, pp. 515--524. Springer Proc. Math. Stat., {\bf 106}
Springer, Tokyo, 2014

\bibitem{Tas15} H. Tasaki, {\it Estimates of antipodal sets in oriented real Grassmann manifolds}, Internat. J. Math. {\bf 26(6)} (2015),  1541008, 12 pp.

\bibitem{Tas17} H. Tasaki, {\it  Sequences of maximal antipodal sets of oriented real Grassmann manifolds II}, Hermitian-Grassmannian submanifolds, pp. 17--26.
Springer Proc. Math. Stat., {\bf 203} (2017).

\bibitem{Wo89} J. A. Wood, {\it Spinor groups and algebraic coding theory}, J. Combin. Theory {\bf 51} (1989), 277--313.

\bibitem{Yu23} J. Yu, {\it Maximal antipodal sets in irreducible compact symmetric spaces},
Transform. Groups {\bf 28(2)} (2023), 987--1000. 

\end{thebibliography}
\end{document}